\documentclass[10pt]{amsart}
\usepackage{amsmath,amssymb,amsrefs}

\numberwithin{equation}{section}

\newtheorem{theorem}[equation]{Theorem}
\newtheorem{lemma}[equation]{Lemma}
\newtheorem{corollary}[equation]{Corollary}
\newtheorem{prop}[equation]{Proposition}

\theoremstyle{definition}
\newtheorem{definition}[equation]{Definition}

\theoremstyle{remark}
\newtheorem{remark}[equation]{Remark}

\begin{document}

\begin{center}
\texttt{Comments, suggestions, corrections, and references are
  most welcomed!}
\end{center}
%\bigskip

\title{Capable two-generator 2-groups of class two}
\author{Arturo Magidin}
\address{Dept. of Mathematical Sciences, The University of
  Montana, Missoula MT 59812}
\email{magidin@member.ams.org}

\subjclass[2000]{Primary 20D15}

\begin{abstract}
A group is called capable if it is a central factor group. We
characterize the capable $2$-generator $2$-groups of class $2$ in
terms of a standard presentation.
\end{abstract}

\maketitle

\section{Introduction}\label{sec:intro}

Recall that a group $G$ is said to be \textbf{capable} if there exists
some group $K$ such that $K/Z(K)$ is isomorphic to~$G$. There are
groups which are not capable (nontrivial cyclic groups being a
well known example), and the question of which finite $p$-groups are
capable plays an important role in their classification. See for
example P.~Hall's comments in his landmark \cite{hallpgroups} (middle
of pp.~137).

In~\cite{baer} Baer characterized the capable groups that are
a direct sum of cyclic groups; the question of which groups are
capable has received renewed attention thanks to results connecting
the question to certain cohomological functors, most notably the
nonabelian tensor square~\cites{beyl,ellis}. There also has been
work studying the restrictions that capability places on the
structure of a group; see for example \cites{heinnikolova,isaacs}. 

More recently, Bacon and Kappe~\cite{baconkappe} characterized the
capable $2$-generator $p$-groups of class two, with $p$ an odd prime,
using the nonabelian tensor square. The author obtained the same
result using other methods~\cite{capable}. The purpose of the present
note is to extend that characterization to the case of $p=2$, yielding
a characterization of capability for
$2$-generated finite groups of class at most two.

\section{A necessary condition}

Notation will be standard; we use the convention that the commutator
of two elements is given by $[x,y]=x^{-1}y^{-1}xy$, and commutators
will be written left-normed, so that $[x,y,z] = [[x,y],z]$, etc. We
will also write $x^y$ to denote $y^{-1}xy$.

The following commutator identities are well known, and may be
verified by direct calculation:

\begin{prop} Let $G$ be any group. Then for all $x,y,z\in G$, $r,s,n\in\mathbb{Z}$,
\begin{itemize}
\item[(a)] $[xy,z] = [x,z]^y[y,z] = [x,z][x,z,y][y,z]$.
\item[(b)] $[x,yz] = [x,z][x,y]^z=[x,z][z,[y,x]][x,y]$.
\item[(c)] $[x^r,y^s] \equiv
  [x,y]^{rs}[x,y,x]^{s\binom{r}{2}}[x,y,y]^{r\binom{s}{2}}\pmod{G_4}$.
\item[(d)] $[y^s,x^r] \equiv
  [x,y]^{-rs}[x,y,x]^{-s\binom{r}{2}}[x,y,y]^{-r\binom{s}{2}}\pmod{G_4}$. 
\item[(e)] $(xy)^n \equiv x^ny^n[y,x]^{\binom{n}{2}}\pmod{G_3}$.
\end{itemize}
Here, $\binom{m}{2} = \frac{m(m-1)}{2}$ for all integers~$m$, and
$G_k$ is the $k$-th term of the lower central series of~$G$.
\label{prop:commident}
\end{prop}

In~\cite{capable}, we established a necessary condition for capability
of a finitely generated $p$-group of class~$k$ in terms of the orders of
the elements of a minimal generating set. In our setting of $2$-groups of class~$2$,
the necessary condition becomes:

\begin{prop} Let $G$ be a group of class two, minimally generated by
  $x_1,\ldots,x_m$, with $x_i$ of order $2^{r_i}$ and satisfying $1\leq
  r_1\leq r_2\leq\cdots\leq r_m$. If $G$ is capable, then $m>1$ and
  $r_m\leq r_{m-1}+1$.
\label{prop:gennecessary}
\end{prop}

When $r_m=r_{m-1}+1$, we can add a condition on the orders of the
commutators $[r_m,r_i]$, $1\leq i\leq m-1$:

\begin{theorem} Let $G$ be a group of class two, minimally generated
  by elements $x_1,\ldots,x_m$ of orders $1< 2^{r_1}\leq\cdots\leq
  2^{r_m}$, respectively, and assume that $m>1$ and
  $r_m=r_{m-1}+1$. If $G$ is capable, then at least one of the
  commutators $[r_m,r_i]$ is of
  order $2^{r_{m-1}}$, for some $i\in\{1,\ldots,m-1\}$. 
\label{th:necessitycommcondition}
\end{theorem}

The theorem will follow from the next lemma:

\begin{lemma} Let $K$ be a nilpotent group of class three, let
  $y_1,\ldots,y_m$ be elements of $K$ which generate $K$ modulo
  $Z(K)$, and assume that $y_i^{2^{r_i}}\in Z(K)$, with
\[ 1 \leq r_1 \leq \cdots \leq r_{m-1}\leq r_m.\]
If there exist integers $0\leq \gamma_i<r_{m-1}$, $i=1,\ldots,m-1$
such that $[y_m,y_i]^{2^{\gamma_i}}$ commutes with $y_i$ and~$y_m$, then
$y_m^{2^{r_{(m-1)}}}\in Z(K)$.
\label{lemma:commcond}
\end{lemma}

\begin{proof} To avoid subscripts of exponents, let
  $\alpha=r_{m-1}$. We want to show that $y_m^{2^{\alpha}}$ commutes
  with $y_i$, $i=1,\ldots,m-1$.
Since $[y_m,y_i]^{2^{\gamma_i}}\in Z(\langle y_i,y_m\rangle)$, it follows that 
\[e  =  [[y_m,y_i]^{2^{\gamma_i}},y_i] =
[y_m,y_i,y_i]^{2^{\gamma_i}}
 =  [[y_m,y_i]^{2^{\gamma_i}},y_m] = [y_m,y_i,y_m]^{2^{\gamma_i}}.\]
Since $\gamma_i < \alpha$, 
$[y_m,y_i,y_i]^{2^{\alpha-1}}$ is a power of
$[y_m,y_i,y_i]^{2^{\gamma_i}}$, and therefore is trivial. Same with
$[y_m,y_i,y_m]^{2^{\alpha-1}}$. We also have
\begin{eqnarray*}
e  &=&  [y_m,y_i^{2^{\alpha}}] =
[y_m,y_i]^{2^{\alpha}}[y_m,y_i,y_i]^{\binom{2^{\alpha}}{2}}\\
 &=&  [y_m,y_i]^{2^{\alpha}}[y_m,y_i,y_i]^{2^{\alpha-1}(2^{\alpha}
  - 1)}
 =  [y_m,y_i]^{2^{\alpha}},
\end{eqnarray*}
so we conclude that $[y_m,y_i]^{2^{\alpha}} = e$ for $i=1,\ldots,m-1$. 
Therefore,
\[\relax[y_m^{2^{\alpha}},y_i]  = 
      [y_m,y_i]^{2^{\alpha}}[y_m,y_i,y_m]^{\binom{2^{\alpha}}{2}}
 =  e.
\]
This proves that $y_m^{2^{\alpha}}\in Z(K)$, as claimed.
\end{proof}

The proof of Theorem~\ref{th:necessitycommcondition} is now immediate:

\begin{proof} Let $K$ be a group such that $G\cong K/Z(K)$. Let 
  $y_1,\ldots,y_m$ be elements of~$K$ that
  project onto $x_1,\ldots,x_m$, respectively. We know that for each
  $i$, the order of $[x_m,x_i]$ divides $2^{r_i}$; if all commutators
  $[x_m,x_i]$ have order strictly less than $2^{r_{(m-1)}}$, then by
  Lemma~\ref{lemma:commcond} we must have $y_m^{r_{(m-1)}}\in Z(K)$, and
  therefore $r_m\leq r_{(m-1)}$. So if $r_m = r_{(m-1)}+1$, then we must
  have that at least one of the commutators $[x_m,x_i]$ is of order at
  least $2^{r_{(m-1)}}$; since this is the highest order it can have,
  we conclude that there exists some $i<m$ with $r_i=
  r_{m-1}$, and with $[x_m,x_i]$ of order exactly $2^{r_{(m-1)}}$, as claimed.
\end{proof}

\section{The classification of two-generated 2-groups of class two}

Since we aim to characterize the capable $2$-generator $2$-groups of
class two, we need a description of these groups.  We modify the
presentation from~\cite{twogrouptensor} to facilitate our own
analysis. 

\begin{theorem}[Theorem~2.5 in~\cite{twogrouptensor}] Let $G$ be a
  finite nonabelian $2$-generator $2$-group of nilpotency class
  two. Then $G$ is isomorphic to exactly one group of the following
  types:
\begin{itemize}
\item[(i)] $\displaystyle G\cong \Bigl\langle a,b\,\Bigm|\,
  a^{2^{\alpha}}\!\!\!=b^{2^{\beta}}\!\!\!=[a,b]^{2^{\gamma}}\!\!\!=[a,b,a]=[a,b,b]=e\Bigr\rangle$, 
where $\alpha$, $\beta$, and $\gamma$ are positive integers satisfying
  $\alpha\geq\beta\geq\gamma$.
\item[(ii)] $\displaystyle G\cong \Bigl\langle a,b\,\Bigm|\,
  a^{2^{\alpha}}\!\!\!=b^{2^{\beta}}\!\!\!=[a,b,a]=[a,b,b]=e,\quad
  a^{2^{\alpha+\sigma-\gamma}}\!\!\!\!\!\!=[a,b]^{2^{\sigma}} \Bigr\rangle$, with
  $\alpha,\beta,\gamma,\sigma$ integers satisfying
  $\beta\geq\gamma>\sigma\geq0$, $\alpha+\sigma\geq 2\gamma$, and $\alpha+\beta+\sigma>3$.
\item[(iii)] $\displaystyle G\cong\Bigl\langle
  a,b\,\Bigm|\,a^{2^{\gamma+1}}\!\!\!\!\!\!=b^{2^{\gamma+1}}\!\!\!\!\!\!%
=[a,b]^{2^{\gamma}}\!\!\!\!=[a,b,a]\!\!=\!\![a,b,b]\!\!=\!\!e,
  a^{2^{\gamma}}\!\!\!\!\!=b^{2^{\gamma}}\!\!\!\!\!=[a,b]^{2^{\gamma-1}}\Bigr\rangle$,
  with $\gamma$ a positive integer.
\end{itemize}
\label{th:classif}
\end{theorem}

\begin{remark} Although in~\cite{twogrouptensor} the commutator
  convention is that $[x,y]=xyx^{-1}y^{-1}$, this difference is
  immaterial because in a nilpotent group of class two we always have
  $xyx^{-1}y^{-1}=x^{-1}y^{-1}xy$.
\end{remark}

The condition $\alpha+\beta+\sigma>3$ in type (ii) is to prevent the
dihedral group of order $8$ from appearing in both types (i) (with
$\alpha=\beta=\gamma=1$) and
(ii) (with $\alpha=2$, $\beta=\gamma=1$, $\sigma=0$).

\begin{remark} In the original statement of Theorem~\ref{th:classif}
  in~\cite{twogrouptensor}, the authors split (ii) into two families,
  one with $\sigma=0$ and one with $\sigma>0$; the former family
  consists of the split metacyclic $2$-groups of class~two.  We could
  even combine types (i) and (ii) into a single expression by using
  the presentation in (ii) and allowing $\gamma=\sigma>0$.
  However, the condition necessary to guarantee
  uniqueness becomes more cumbersome, and since we will deal
  with groups of type (i) separately below anyway, we have kept
  the distinction.
\end{remark}

I will refer to the three types informally as follows: groups of type (i) will be the
\textit{``coproduct type''} groups, since they arise as the coproduct
of two cyclic groups in the subvariety of all nilpotent groups of
class two in which $2^{\gamma}$-powers are central. Groups of type (ii) will be the
\textit{``general type''}
groups, since they give rise to all groups between split metacyclic
and coproducts. When $p$ is odd, the two-generator $p$-groups of class
two are counterparts to these two families; the groups of type~(iii)
above have no such counterparts for odd primes, so groups of
type~(iii) will be called of \textit{``exceptional type.''} 

\section{The exceptional type case}

The exceptional type is straightforward: a group of
exceptional type is never capable.

\begin{theorem} If $G$ is of exceptional type (that is, presented as in
  Theorem~\ref{th:classif}(iii)), then $G$ is not capable.
\end{theorem}

\begin{proof} Let $K$ be a group of class three, and assume that
  $K/N\cong G$, with $N$ a central subgroup. We want to show that
 $N\neq Z(K)$. 
 Let $x,y\in K$ project down to $a$ and~$b$ in~$G$,
  respectively. Since $a^{2^{\gamma}}=b^{2^{\gamma}}$ in $G$, we must
  have $x^{2^{\gamma}}y^{-2^{\gamma}}\in N\subset Z(K)$. Since
 $y^{2^{\gamma}}$ commutes with $y$, so does the product
 $(x^{2^{\gamma}}y^{-2^{\gamma}})y^{2^{\gamma}}=x^{2^{\gamma}}$. As
 $x^{2^{\gamma}}$ also commutes with $x$, we conclude that it is
 in fact central in~$K$. Since $a^{2^{\gamma}}\neq e$, we must have $x^{2^{\gamma}}\not\in
 N$. This proves $N\neq Z(K)$, as desired.
\end{proof}

\section{Normal forms and the coproduct type case}

Our method consists of constructing explicit witnesses to the
capability of the given groups when appropriate. As we attempt the
construction, the cases where the group is not capable will become
evident by the appearance of undesired relations in the potential witness. In
essence, we are constructing a simplified version of the ``generalised
extension of G;'' see Theorem III.3.9 in~\cite{beyltappe}. Our
starting point will be the 3-nilpotent product of cyclic groups (see
\cites{struikone,struiktwo,capable} for details).  We restrict the general
definition to the case we are interested in:

\begin{definition} Let $A_1,\ldots,A_t$ be cyclic groups. The
  $3$-nilpotent product of the $A_i$, denoted
$A_1 \amalg^{{\germ N}_3}\cdots\amalg^{{\germ N}_3} A_t$,
is defined to be $F/F_{4}$, where $F$ is the free product of the $A_i$,
$F=A_1*\cdots*A_t$, and $F_{4}$ is the fourth term of the lower central
series of~$F$. 
\end{definition}

We will mostly be concerned with the case $t=2$, so we restrict our
presentation below to that situation.
The normal form given in~\cite{struikone} for the $3$-nilpotent
product of cyclic $2$-groups seems to be difficult to use in most of
our cases. That choice of normal form was made to facilitate the
description of the multiplication table, and we will not need the
multiplication table. We will therefore also provide an alternative normal
form which we will use to study some of the cases below. 

\begin{theorem}[Struik, Theorem~4 in~\cite{struikone}]
Let $a$ and $b$ generate cyclic groups of order $2^{\alpha}$ and
$2^{\beta}$, respectively, with $\alpha\geq\beta\geq 1$. Let
$G=\langle a\rangle\amalg^{{\germ N}_3}\langle b\rangle$ be their
$3$-nilpotent product. Then every element of $G$ may be expressed
uniquely in the form
\begin{equation}
a^r b^s [a,b]^t [a^2,b]^u [a,b^2]^v
\label{eq:struiknormalform}
\end{equation}
where $r$ is unique modulo $2^{\alpha}$; $s$ is unique modulo
$2^{\beta}$; $t$ is unique modulo $2^{\beta+1}$; $v$ is unique modulo
$2^{\beta-1}$; and $u$ is unique modulo $2^{\beta}$ if
$\alpha\neq\beta$, and unique modulo $2^{\beta-1}$ if $\alpha=\beta$.
\label{th:struiknormalform}
\end{theorem}

\begin{corollary}[cf. Theorems 5.1 and 5.2 in~\cite{capable}]
Let $G$ be as above. The center of $G$ is generated by
$a^{2^{\beta+1}}$, $[a,b]^2[a^2,b]^{-1}$, and
$[a,b]^2[a,b^2]^{-1}$. Therefore,
\[ G/Z(G) \cong \bigl\langle x,y \,\bigm|\,
x^{2^{\max\{\alpha,\beta+1\}}}=y^{2^{\beta}}=[x,y]^{2^{\beta}}=[x,y,x]=[x,y,y]=e\bigr\rangle.\]
\label{cor:centralquotientone}
\end{corollary}

Suppose, however, that want to make
$[a,b]^{2^{\gamma}}$ is central, for some $\gamma<\beta$, so the
central quotient will have $[a,b]$ of order $2^{\gamma}$. In that
case, it is more convenient to switch to a normal form that uses basic
commutators on $a$ and~$b$. We have the following:

\begin{theorem}Let $\langle a\rangle$, $\langle b\rangle$ be cyclic
  groups of order $2^{\alpha}$ and $2^{\beta}$, respectively, with
  $\alpha\geq\beta\geq1$. Let $G = \langle a\rangle\amalg^{{\germ
  N}_3}\langle b\rangle$, and let $\gamma$ be an integer, $1\leq
  \gamma < \beta$. Let $N$ be the central subgroup of $G$ generated by
$[a,b,a]^{2^{\gamma}}  =  [a,b]^{-2^{\gamma+1}}[a^2,b]^{2^{\gamma}}$
  and 
$[a,b,b]^{2^{\gamma}}  =  [a,b]^{-2^{\gamma+1}}[a,b^2]^{2^{\gamma}}$.
Then every element of $K=G/N$ can be written uniquely as 
$k=a^r b^s [a,b]^t[a,b,a]^u[a,b,b]^v$
(abusing notation and writing $a$ instead of $aN$, etc.),
where $r$ is unique modulo $2^{\alpha}$, $s$
and~$t$ are unique modulo $2^{\beta}$, and $u$ and $v$ are unique
modulo $2^{\gamma}$. 
\label{th:mynormalform}
\end{theorem}

\begin{proof} We may rewrite any element of $G$ as
\begin{equation}
g = a^r b^s [a,b]^t [a,b,a]^u [a,b,b]^v
\label{eq:preliminarynormal}
\end{equation}
by using the normal form from Theorem~\ref{th:struiknormalform},
the identities $[a^2,b]=[a,b]^2[a,b,a]$ and $[a,b^2]=[a,b]^2[a,b,b]$,
and the fact that $[G,G]$ is abelian. Note that:
\[ e =
   [a,b^{2^{\beta}}]=[a,b]^{2^{\beta}}[a,b,b]^{\binom{2^{\beta}}{2}} =
   [a,b]^{2^{\beta}}[a,b,b]^{2^{\beta-1}};\]
the last equality since $2\beta-1\geq \beta$, and
   $[a,b,b]^{2^{\beta-1}}$ is of order $2$. 
A straightforward calculation shows that we may choose the exponents
in $(\ref{eq:preliminarynormal})$ so that $r$ is unique modulo
$2^{\alpha}$; $s$ is unique modulo $2^{\beta}$, $t$ is unique modulo
$2^{\beta+1}$; $u$ is either unique modulo $2^{\beta}$ if
$\alpha\neq\beta$, or satisfies $0\leq u < 2^{\beta-1}$ if
$\alpha=\beta$; and $v$ satisfies $0\leq v < 2^{\beta-1}$. Now simply
   note that $[a,b]^{2^{\beta}}\in N$,
   $[a,b]^{2^{\beta}}=[a,b,b]^{2^{\beta-1}}$ (and
   $[a,b]^{2^{\beta}}=[a,b,a]^{2^{\beta-1}}$ if $\alpha=\beta$), 
and we obtain the normal form
   described by taking the quotient.
\end{proof}

\begin{corollary} Let $K$ be as in
   Theorem~\ref{th:mynormalform}. Then $Z(K)$ is generated by
   $a^{2^{\beta}}$, $[a,b]^{2^{\gamma}}$, $[a,b,a]$, and $[a,b,b]$. Therefore,
\[ K/Z(K) \cong \Bigl\langle x,y\,\Bigm|\, x^{2^{\beta}} =
   y^{2^{\beta}}=[x,y]^{2^{\gamma}}=[x,y,x]=[x,y,y]=e\Bigr\rangle.\]
%% Corrected the size of the \rangle
\label{cor:centralquotienttwo}
\end{corollary}

\begin{proof} Using the normal form and
   Proposition~\ref{prop:commident}, it is straightforward to verify
   that the elements given generate the center; the description of the
   quotient follows by mapping $a$ to $x$ and $b$ to $y$.
\end{proof}

These corollaries, together with Prop.~\ref{prop:gennecessary} and
Theorem~\ref{th:necessitycommcondition}, yield the coproduct type case:

\begin{corollary} Let $G$ be a two-generated $2$-group of class two
  and of coproduct type; that is, 
  presented as in Theorem~\ref{th:classif}(i). Then $G$ is capable if
  and only if $\alpha=\beta$, or $\alpha=\beta+1$ and $\gamma=\beta$.
\end{corollary}

\section{General type, $\gamma<\beta$}

We turn now to the case of groups of general type that have
$\gamma<\beta$; by Theorem~\ref{th:necessitycommcondition}, 
% CORRECTION: from ``this implies that $\alpha=\beta$.'' to
we may assume that $\alpha=\beta$. 

Let $K$ be the group
described in Theorem~\ref{th:mynormalform}. Let $N$ be the subgroup
of~$K$ generated by
$[a^{2^{\alpha+\sigma-\gamma}}[a,b]^{-2^{\sigma}},a]$ and 
$[a^{2^{\alpha+\sigma-\gamma}}[a,b]^{-2^{\sigma}},b]$.
We will show that $K/N$ has $G$ as a central quotient, provided
    $\alpha>\gamma+1$. 
First we give a better description of $N$. We rewrite the elements
above in normal form:
\begin{eqnarray*}
\relax [a^{2^{\alpha+\sigma-\gamma}}[a,b]^{-2^{\sigma}},a] & = &
       [[a,b]^{-2^{\sigma}},a] = [a,b,a]^{-2^{\sigma}}.\\
\relax [a^{2^{\alpha+\sigma-\gamma}}[a,b]^{-2^{\sigma}},b] & = &
       [a^{2^{\alpha+\sigma-\gamma}},b][[a,b]^{-2^{\sigma}},b]\\
& = &
       [a,b]^{2^{\alpha+\sigma-\gamma}}[a,b,a]^{\binom{2^{\alpha+\sigma-\gamma}}{2}}[a,b,b]^{-2^{\sigma}}\\
& = &
       [a,b]^{2^{\alpha+\sigma-\gamma}}[a,b,a]^{(2^{\alpha+\sigma-\gamma-1})(2^{\alpha+\sigma-\gamma}-1)}
[a,b,b]^{-2^{\sigma}}.\\
\end{eqnarray*}
Since $\alpha+\sigma\geq 2\gamma$, we must have
$\alpha+\sigma-\gamma-1 \geq \gamma-1\geq \sigma$. Thus, the subgroup
$N$ is generated by
$[a,b]^{2^{\alpha+\sigma-\gamma}}[a,b,b]^{-2^{\sigma}}$ and $[a,b,a]^{2^{\sigma}}$.
Since both elements are central, we have that an element written in
   normal form $a^r b^s
[a,b]^t [a,b,a]^u [a,b,b]^v \in K$ will lie in $N$ if and only if
$r\equiv s\equiv 0 \pmod{2^{\alpha}}$, $u\equiv v\equiv 0
\pmod{2^{\sigma}}$, and $t+2^{\alpha-\gamma}v\equiv
0\pmod{2^{\alpha}}$. Note that the last expression is well defined,
since $v$ is well defined modulo $2^{\gamma}$.

If $G$ is to be the central quotient of $K/N$, then we must have that
$kN\in Z(K/N)$ if and only if $k$ lies in $\langle
a^{2^{\alpha+\sigma-\gamma}}[a,b]^{-2^{\sigma}}\rangle Z(K)$.
Let $k\in K$ be an arbitrary element, written in normal form as
$k = a^r b^s [a,b]^t [a,b,a]^u [a,b,b]^v$, and assume that
$kN$ is central in $K/N$. That means that both $[k,a]$ and $[k,b]$
lie in $N$. We have:
\begin{eqnarray}
\relax [k,a] & = & [a,b]^{-s}[a,b,a]^t
       [a,b,b]^{-\binom{s}{2}}.\label{eq:kbracketa}\\
\relax [k,b] & = & [a,b]^r
       [a,b,a]^{\binom{r}{2}}[a,b,b]^{rs+t}.\label{eq:kbracketb}
\end{eqnarray}
Therefore, we must have:
\begin{eqnarray*}
-s - 2^{\alpha-\gamma}\binom{s}{2} = -s\left(1+2^{\alpha-\gamma-1}(s-1)\right)&\equiv& 0 \pmod{2^{\alpha}}\\
r + 2^{\alpha-\gamma}(rs+t) &\equiv& 0 \pmod{2^{\alpha}}\\
t\equiv \binom{s}{2}\equiv \binom{r}{2} \equiv rs+t &\equiv & 0
\pmod{2^{\sigma}}.
\end{eqnarray*}
From the first two conditions we obtain that $r$ and $s$ are even. If
$\alpha>\gamma+1$, then the first congruence gives $s\equiv
0\pmod{2^{\alpha}}$, and the second reduces to
$r+2^{\alpha-\gamma}t\equiv 0 \pmod{2^{\alpha}}$. Since $t$ is
divisible by $2^{\sigma}$, we have $t=m2^{\sigma}$ and
$r\equiv-m2^{\alpha+\sigma-\gamma}\pmod{2^{\alpha}}$ for some
integer~$m$. The remaining congruences now follow. Thus:
\[ k = \left( a^{2^{\alpha+\sigma-\gamma}}[a,b]^{2^\sigma}\right)^{-m}
   [a,b,a]^{u+\binom{-m}{2}2^{\alpha+2\sigma-\gamma}}[a,b,b]^v,\]
so $k$ lies in the desired subgroup. 

If $\alpha=\gamma+1$, then from $\alpha+\sigma\geq\gamma>\sigma$ we
conclude that $\sigma+1=\gamma$; as far as our congruences are
concerned, $r=t=0$ and $s=2^{\alpha-1}=2^{\gamma}$ is also a
solution but does not fit into the subgroup we want. This 
suggests the following observation:

\begin{lemma} Let $K$ be a nilpotent group of class~$3$, and let $x$
  and $y$ be elements of~$K$. Assume that
  $\alpha>1$ is an integer such that $x^{2^{\alpha}}$,
  $[x,y]^{2^{\alpha-1}}$, and $x^{2^{\alpha-1}}[x,y]^{-2^{\alpha-2}}$
  centralize $\langle x,y\rangle$.
 Then $y^{2^{\alpha-1}}$ commutes with $x$.
\end{lemma}

\begin{proof}
Since $[x,y]^{2^{\alpha-1}}\in Z(\langle x,y\rangle)$, we must have
$e = [x,y,x]^{2^{\alpha-1}}\!\!\!\! = [x,y,y]^{2^{\alpha-1}}$.
Therefore, $[x,y,y]^{-2^{\alpha-2}} = [x,y,y]^{2^{\alpha-2}}$.
In addition, since $[x,y]^{2^{\alpha-1}}$ centralizes $\langle
x,y\rangle$, we have
$[x,y]^{-2^{\alpha-2}}\equiv [x,y]^{2^{\alpha-2}}\pmod{Z(\langle x,y\rangle)}$. 

From the other centralizing elements we deduce that:
\[ e = [x^{2^{\alpha}},y] =
   [x,y]^{2^{\alpha}}[x,y,x]^{\binom{2^{\alpha}}{2}} =
   [x,y]^{2^{\alpha}}.\]
and also that
$[x^{2^{\alpha-1}}[x,y]^{2^{\alpha-2}},x] = [x,y,x]^{2^{\alpha-2}}\!\!\!\! = e$.
In addition,
\begin{eqnarray*}
e  &=& [x^{2^{\alpha-1}}[x,y]^{2^{\alpha-2}},y]
 =
[x,y]^{2^{\alpha-1}}[x,y,x]^{\binom{2^{\alpha-1}}{2}}[x,y,y]^{2^{\alpha-2}}\\
& = & [x,y]^{2^{\alpha-1}}[x,y,y]^{2^{\alpha-2}} = [x,y]^{2^{\alpha-1}}[x,y,y]^{-2^{\alpha-2}}.
\end{eqnarray*}
From these equations we obtain:
\begin{eqnarray*}
\relax[x,y^{2^{\alpha-1}}] & = &
      [x,y]^{2^{\alpha-1}}[x,y,y]^{\binom{2^{\alpha-1}}{2}}
 =  [x,y]^{2^{\alpha-1}}[x,y,y]^{2^{\alpha-2}(2^{\alpha-1}-1)}\\
& = &
      [x,y]^{2^{\alpha-1}}[x,y,y]^{-2^{\alpha-2}}[x,y,y]^{2^{2\alpha-3}}
 =  e.
\end{eqnarray*}
The last equality holds since $\alpha>1$, so $[x,y,y]^{2^{2\alpha-3}}=e$.
\end{proof}

Thus we obtain:

\begin{corollary}Let $G$ be a group of general type (that is, presented
  as in Theorem~\ref{th:classif}(ii)) with $\gamma<\beta$. Then $G$ is capable if and only if
  $\alpha=\beta$ and $\gamma<\beta-1$. 
\end{corollary}

%% Fixed proof below; replaced some K's with H's.
\begin{proof} We have seen the necessity of $\alpha=\beta$. 
If $\gamma+1<\beta=\alpha$, we have shown that the group $K/N$
  described above satisfies $(K/N)/Z(K/N)\cong G$, so $G$ is
  capable. Assume then that $\gamma+1=\beta=\alpha$ (and so
  $\sigma=\gamma-1=\alpha-2$), and let $H$ be a group such that $H/N\cong G$,
  with $N\subset Z(H)$. We want to show that $N\neq Z(H)$. Note that
  $H$ must be of class at most three. Let $x,y\in H$ project down
  to $a$ and $b$, respectively. In $H$ we have that $x^{2^{\alpha}}$,
  $[x,y]^{2^{\gamma}}=[x,y]^{2^{\alpha-1}}$, and
  $x^{2^{\alpha+\sigma-\gamma}}[x,y]^{-2^{\sigma}} =
  x^{2^{\alpha-1}}[x,y]^{-2^{\alpha-2}}$ lie in $N$, and therefore are
  central. From the lemma we conclude that $y^{2^{\alpha-1}}\in
  Z(H)$. Since $b^{2^{\alpha-1}}\neq e$, we must have $N\neq Z(H)$,
  which proves $G$ cannot be capable.
\end{proof}

\section{General type, $\gamma=\beta$}

Although the description of groups of general type admits the case
with $\alpha<\beta$, such a group cannot be capable: if
$\alpha<\beta$, then Theorem~\ref{th:necessitycommcondition} implies
that $\gamma=\alpha$, and this makes $\alpha+\sigma\geq 2\gamma$
impossible. So we must have $\alpha\geq\beta$. And if $\gamma=\beta$,
then we must have $\alpha>\beta$, for otherwise we again cannot
satisfy the condition $\alpha+\sigma\geq 2\gamma$. Thus, the case we
are to consider now is $\alpha=\beta+1=\gamma+1$.  This in turns
implies $\sigma=\alpha-2$.  In this situation,
$\alpha+\sigma-\gamma=\beta$.

Since $\gamma=\beta$, it is preferable to use the normal form in
Theorem~\ref{th:struiknormalform}. Let $a$ generate a cyclic group of
order $2^{\alpha}=2^{\beta+1}$, and let $b$ generate a cyclic group of
order $2^{\beta}$. Let $K = \langle a\rangle \amalg^{{\germ N}_3}
\langle b\rangle$, and
let $N$ be the subgroup of $K$ generated by the elements
$[a^{2^{\beta}}[a,b]^{-2^{\beta-1}},a]$ and $[a^{2^{\beta}}[a,b]^{-2^{\beta-1}},b]$.
We want to show that if $k\in K$ has $[k,a],[k,b]\in N$, then $k$ lies
in the subgroup $\langle a^{2^{\beta}}[a,b]^{-2^{\beta-1}}\rangle Z(K)$, and so
deduce that $G$ is the central quotient of $K/N$.

First, we calculate the normal forms of the generators of~$N$:
\begin{eqnarray*}
\relax[a^{2^{\beta}}[a,b]^{-2^{\beta-1}},a] & = &
      [a,b,a]^{-2^{\beta-1}}
 =  \left([a,b]^{-2}[a^2,b]\right)^{-2^{\beta-1}}\\
& = & [a,b]^{2^{\beta}}[a^2,b]^{-2^{\beta-1}}
 =  [a,b]^{\pm2^{\beta}}[a^2,b]^{\pm2^{\beta-1}}.\\
\relax [a^{2^{\beta}}[a,b]^{-2^{\beta-1}},b] & = &
      [a,b]^{2^{\beta}}[a,b,a]^{\binom{2^{\beta}}{2}}[a,b,b]^{-2^{\beta-1}}\\
& = &
      [a,b]^{2^{\beta}}\left([a,b]^{-2}[a^2,b]\right)^{\binom{2^{\beta}}{2}}
      \left([a,b]^{-2}[a,b^2]\right)^{-2^{\beta-1}}\\
& = & [a,b]^{2^{\beta}-2\binom{2^{\beta}}{2} +
      2^{\beta}}[a^2,b]^{\binom{2^{\beta}}{2}}
      [a,b^2]^{-2^{\beta-1}}\\
& = & [a,b]^{-2^{\beta}(2^{\beta}-1)}[a^2,b]^{2^{\beta-1}(2^{\beta}-1)}\\
& = & [a,b]^{2^{\beta}}[a^2,b]^{-2^{\beta-1}}
 =  [a,b]^{\pm2^{\beta}}[a^2,b]^{\pm2^{\beta-1}}.
\end{eqnarray*}
That is, $N$ is central and cyclic of order $2$, generated by
$[a,b]^{\pm2^{\beta}}[a^2,b]^{\pm2^{\beta-1}}$. The liberty in signs
is because both $[a,b]^{2^{\beta}}$ and $[a^2,b]^{2^{\beta-1}}$ are of
order two.

An arbitrary element $a^r b^s [a,b]^t [a^2,b]^u [a,b^2]^v$ will
be in $N$ if and only if $r\equiv 0 \pmod{2^{\beta+1}}$, $s\equiv
0\pmod{2^{\beta}}$, $u\equiv v\equiv 0\pmod{2^{\beta-1}}$ (note that
$[a,b^2]$ is of order $2^{\beta-1}$, while $[a^2,b]$ is of order
$2^{\beta}$), and $t + 2u\equiv 0 \pmod{2^{\beta+1}}$. This last
expression is well defined, since $u$ is well defined modulo
$2^{\beta}$.

Now let $k\in K$ be given by $k=a^r b^s [a,b]^t [a^2,b]^u [a,b^2]^v$,
and assume that both commutators $[k,a]$ and $[k,b]$ lie in~$N$.  We
may start from $(\ref{eq:kbracketa})$ and~$(\ref{eq:kbracketb})$, and
substitute the values of $[a,b,a]$ and $[a,b,b]$. We obtain:
\begin{eqnarray*}
\relax[k,a] & = & [a,b]^{2\binom{s}{2}-s-2t}[a^2,b]^t[a,b^2]^{-\binom{s}{2}},\\
\relax[k,b] & = & [a,b]^{r-2\binom{r}{2}-2(rs+t)}[a^2,b]^{\binom{r}{2}}[a,b^2]^{rs+t}.
\end{eqnarray*}
Since both $[k,a]$ and $[k,b]$ lie in~$N$, we must have $r$ and $s$
even; in addition, we have 
$ \binom{s}{2}\equiv \binom{r}{2}\equiv 0
\pmod{2^{\beta-1}}$.
We conclude that
$r\equiv s\equiv 0 \pmod{2^\beta}$. Since $s$ is only defined modulo
$2^{\beta}$, we may take $s=0$.

Finally, we must also have 
\[ r-2\binom{r}{2}-2t + 2\binom{r}{2} \equiv r-2t \equiv 0
\pmod{2^{\beta+1}}.\] Since $r\equiv 0 \pmod{2^{\beta}}$ and $t\equiv
0 \pmod{2^{\beta-1}}$, we have the following possibilities for $r$ and
$t$: $r=t=0$; or $r=0$, $t=2^{\beta}$; or $r=2^{\beta}$, $t=\pm
2^{\beta-1}$; or $r=2^{\beta}$, $t=2^{\beta}\pm 2^{\beta-1}$.

If $r=0$, we obtain $k\in K_3\subset Z(K)$. If $r=2^{\beta}$, then
we have $k$ in the coset of $a^{2^{\beta}}[a,b]^{-2^{\beta-1}}$ modulo
$Z(K)$. In any case, we conclude that $kN$ is central in $K/N$,
if and only if $k\in\langle
a^{2^{\beta}}[a,b]^{-2^{\beta-1}}\rangle Z(K)$.  So $G$ is the central
quotient of $K/N$. 

\begin{corollary}
Let $G$ be a group of general type (that is, presented
  as in Theorem~\ref{th:classif}(ii)) with $\gamma=\beta$. Then $G$ is
  capable if and only if $\alpha=\beta+1$ and $\sigma=\beta-1$.
\end{corollary}

\section{Conclusion}

We summarize our results in the following theorem:

\begin{theorem}
Let $G$ be a $2$-generated $2$-group of class two, presented as in
Theorem~\ref{th:classif}.  Then $G$ is capable if and only if one of
the following hold:
\begin{itemize}
\item[(a)] $G$ is of type (i) and $\alpha=\beta$; or
\item[(b)] $G$ is of type (i), and $\alpha=\beta+1=\gamma+1$; or
\item[(c)] $G$ is of type (ii), $\alpha=\beta$, and $\gamma<\beta-1$; or
\item[(d)] $G$ is of type (ii), and $\alpha=\beta+1=\gamma+1=\sigma+2$.
\end{itemize}
In particular, if $G$ is of type (iii) (exceptional type),
then $G$ is not capable.
\label{th:capablecharacterization}
\end{theorem}


\section*{References}

\begin{biblist}
%\bibselect{bibliog}
\bib{baconkappe}{article}{
  author={Bacon, Michael~R.},
  author={Kappe, Luise-Charlotte},
  title={On capable $p$-groups of nilpotency class two},
  date={2003},
  journal={Illinois J. Math.},
  number={1/2},
  volume={47},
  pages={49\ndash 62},
}
\bib{baer}{article}{
  author={Baer, Reinhold},
  title={Groups with preassigned central and central quotient group},
  date={1938},
  journal={Trans. Amer. Math. Soc.},
  volume={44},
  pages={387\ndash 412},
}
\bib{beyltappe}{book}{
  author={Beyl, F.~Rudolf},
  author={Tappe, J\"urgen},
  title={Group extensions, representations, and the Schur multiplicator},
  date={1982},
  series={Lecture Notes in Mathematics},
  volume={958},
  publisher={Springer-Verlag},
  review={\MR {84f:20002}},
}
\bib{beyl}{article}{
  author={Beyl, F.~Rudolf},
  author={Felgner, Ulrich},
  author={Schmid, Peter},
  title={On groups occurring as central factor groups},
  date={1979},
  journal={J. Algebra},
  volume={61},
  pages={161\ndash 177},
  review={\MR {81i:20034}},
}
\bib{ellis}{article}{
  author={Ellis, Graham},
  title={On the capability of groups},
  date={1998},
  journal={Proc. Edinburgh Math. Soc.},
  number={41},
  pages={487\ndash 495},
  review={\MR {2000e:20053}},
}
\bib{hallpgroups}{article}{
  author={Hall, P.},
  title={The classification of prime-power groups},
  date={1940},
  journal={J. Reine Angew. Math},
  volume={182},
  pages={130\ndash 141},
  review={\MR {2,211b}},
}
\bib{heinnikolova}{article}{
  author={Heineken, Hermann},
  author={Nikolova, Daniela},
  title={Class two nilpotent capable groups},
  date={1996},
  journal={Bull. Austral. Math. Soc.},
  volume={54},
  number={2},
  pages={347\ndash 352},
  review={\MR {97m:20043}},
}
\bib{isaacs}{article}{
  author={Isaacs, I. M.},
  title={Derived subgroups and centers of capable groups},
  date={2001},
  journal={Proc. Amer. Math. Soc.},
  volume={129},
  number={10},
  pages={2853\ndash 2859},
  review={\MR {2002c:20035}},
}
\bib{twogrouptensor}{article}{
  author={Kappe, Luise-Charlotte},
  author={Visscher, Matthew P.},
  author={Sarmin, Nor Haniza},
  title={Two-generator two-groups of class two and their nonabelian tensor square},
  date={1999},
  journal={Glasgow Math. J.},
  volume={41},
  pages={417\ndash 430},
  review={\MR {2000k:20037}},
}
\bib{capable}{article}{
  author={Magidin, Arturo},
  title={Capability of nilpotent products of cyclic groups},
  eprint={{arXiv:math.GR/0403188}},
  note={Submitted},
}
\bib{struikone}{article}{
  author={Struik, Ruth~Rebekka},
  title={On nilpotent products of cyclic groups},
  date={1960},
  journal={Canad. J. Math.},
  volume={12},
  pages={447\ndash 462},
  review={\MR {22:\#11028}},
}
\bib{struiktwo}{article}{
  author={Struik, Ruth~Rebekka},
  title={On nilpotent products of cyclic groups II},
  date={1961},
  journal={Canad. J. Math.},
  volume={13},
  pages={557\ndash 568},
  review={\MR {26:\#2486}},
}

\end{biblist}
\end{document}